\documentstyle[12pt]{article}
\begin{document}
\def\cross{{>\!\!\!\triangleleft}}

\vspace*{.5cm}

\begin{center}
{\bf ON THE ALGEBRAIC STRUCTURE OF DIFFERENTIAL CALCULUS \\
ON QUANTUM GROUPS}
\end{center}
\begin{center}
{\sc O.\,V.\,Radko and A.\,A.\,Vladimirov}
\end{center}
\begin{center}
{\em Bogoliubov Laboratory of Theoretical Physics, \\
Joint Institute for Nuclear Research, \\
Dubna, Moscow region 141980, Russia}
\end{center}

\vspace{.5cm}

{\small 
Intrinsic Hopf algebra structure of the Woronowicz differential
complex is shown to generate quite naturally a bicovariant algebra
of four basic objects within a differential calculus on quantum groups 
-- coordinate functions, differential 1-forms, Lie derivatives, and 
inner derivations -- as the cross-product algebra of two mutually dual
graded Hopf algebras. This construction, properly taking into account 
Hopf-algebraic properties of Woronowicz's bicovariant calculus, provides a 
direct proof of the Cartan identity and of many other useful relations. A 
detailed comparison with other approaches is also given.}

\vspace{1cm}

{\bf I. \ Introduction}

\vspace{.1cm}

Non-commutative differential calculus on quantum groups initiated
and thoroughly worked out by Woronowicz~\cite{Wo} is up to now a subject of
active discussions and development. Though meeting some 
problems~\cite{CWSWW,SWZ1,IP1,IP2,AIP} 
with non-classical dimensionalities of spaces of higher-order differential 
forms (which, in its turn, stimulated very interes\-ting alternative 
approaches~\cite{SWZ2,IP2,FP}), original Woronowicz's construction remains 
highly attractive due to both its rich algebraic structure and useful 
applications. Probably, the best known realization of this scheme is 
bicovariant differential calculus on the $GL_q(N)$ quantum 
groups~\cite{Ju,SWZ1,AC}.

Closely related but somewhat parallel to Woronowicz's construction is 
another project~\cite{SWZ1,SWZ2,SWZ3,Is} that, in particular, has produced 
a bicovariant algebra of four types of elements: functions on a quantum 
group, differential forms, Lie derivatives along vector fields, and inner 
derivations -- by exact analogy with classical differential geometry. 
However, this scheme, as it is, does not seem to be fully motivated by the 
Hopf-algebraic nature of non-commutative differential calculus.

In the present paper, we suggest an extension of Woronowicz's axiomatics
which naturally involves Lie derivatives and inner derivations in a way 
that respects the Hopf algebra structure of the whole scheme. Actually,
in the framework of Woronowicz's noncommutative differential calculus
~\cite{Wo,AC} one deals with the differential complex
\begin{equation}
A\stackrel{d}{\longrightarrow}\Gamma\stackrel{d}{\longrightarrow}
\Gamma^2\longrightarrow\ldots \ \ , 
\label{1}
\end{equation}
where $A$ is a Hopf algebra (\,of functions on a quantum group\,), 
 $\Gamma$ is its bicovariant bimodule,
$\Gamma^2\equiv \Gamma\wedge\Gamma$ is its second wedge power,
and so on. Exterior differential map $d: \Gamma^n
\rightarrow\Gamma^{n+1}$ is assumed to obey the Leibniz rule 
\begin{equation}
d(ab)=(da)b+a\,db 
\label{2}
\end{equation}
and the nilpotency condition $d\circ d = 0 $.   
Brzezinski~\cite{Br} has shown that 
\begin{equation}
\Gamma^{\wedge}\doteq A\oplus\Gamma\oplus\Gamma^2\oplus\ldots 
\label{3}
\end{equation}
also becomes a (graded) Hopf algebra with respect to (wedge)
multiplication and na\-tu\-ral definitions of coproduct and antipode.
In what follows, we want to demonstrate how this Hopf structure
can be used to build an associative noncommutative 
bicovariant algebra containing functions, differential
forms, Lie derivatives and inner derivations. Similar algebras
have been introduced and studied by several 
authors~\cite{SWZ1,SWZ3,Sc1,Is} (and the idea to a use cross-product for
constructing bicovariant differential calculus is due to~\cite{SWZ2}).  
Probably, the closest to ours is 
the approach by P.\,Schupp~\cite{Sc1}. However, some of our 
results and, especially, starting points appear to be different. So, we 
propose the construction described below as entirely
Hopf-algebra motivated (and, we believe, natural) new approach to the
problem.

\vspace{.3cm}

{\bf II. Cross-product of dual Hopf algebras}

\vspace{.1cm}

Notions of mutually dual Hopf algebras and their cross-product will
 actively be used throughout this paper. Let us recall the corresponding
terminology and basic definitions~\cite{Sw,Ab,CP,Ma1,SWZ2}. 
Let  $A$ be a Hopf algebra with associative multiplication, coassociative
coproduct
\begin{equation}
\Delta : A\longrightarrow A\otimes A\,, \ \ \
\Delta (a)\doteq a_{(1)}\otimes a_{(2)}\,, \ \ \
\Delta(ab)=\Delta(a)\Delta(b)
\label{4}
\end{equation}
(we will use the notation
\begin{equation}
a_{(1)}\otimes a_{(2)}\otimes a_{(3)}
\doteq(\Delta\otimes id)\circ\Delta(a)
=(id\otimes \Delta)\circ\Delta(a)\,,
\label{5}
\end{equation}
and so on, for multiple coproducts), a counit
\begin{equation}
\varepsilon : A\longrightarrow C\,, \ \ \
\varepsilon (ab)=\varepsilon (a)\varepsilon (b)\,, \ \ \
\varepsilon (a_{(1)})\,a_{(2)}=a_{(1)}\varepsilon(a_{(2)})=a\,,
\label{6}
\end{equation}
and an invertible antipode
$$
S : A\longrightarrow A\,, \ \ \ S(ab)=S(b)S(a)\,, \ \ \
\Delta(S(a))=S(a_{(2)})\otimes S(a_{(1)})\,,
$$
\begin{equation}
\varepsilon (S(a))=\varepsilon (a)\,, \ \ \
S(a_{(1)})a_{(2)}=a_{(1)}S(a_{(2)})=\varepsilon (a)\,.
\label{7}
\end{equation}
Algebra $A^{*}$ is a Hopf dual of $A$ with  
   $<\! \cdot , \cdot \!> : A^{*}\otimes A\longrightarrow C $ 
being a duality map, if
$$
<\! xy\,, a \!> = <\! x\otimes y, \Delta (a)\!>\,, \ \ \
<\! x\,, ab\!>=<\! \Delta(x), a\otimes b \!>\,,
$$
\begin{equation}
\label{8}
<\!x, 1>=\varepsilon (x)\,, \ \ \
<\! 1\,, a \!>= \varepsilon (a)\,, \ \ \
<\!S(x),a\!>=<\! x,S(a) \!>\,.
\end{equation}
Here and below $a,b\in A\,, x,y \in A^{*}$.

One can define left and right covariant actions $A^*\triangleright A$ and
$A\triangleleft A^*$ by
\begin{equation}
x\triangleright a= a_{(1)}<\!x,a_{(2)}\!>\,, \ \ \  
a\triangleleft x= a_{(2)}<\!x,a_{(1)}\!>\,. \ \ \ 
\label{9}
\end{equation}
As usual, left and right actions imply
\begin{equation}
\label{10}
(xy)\triangleright a=x\triangleright(y\triangleright a)\,, \ \ \
a\triangleleft(xy)=(a\triangleleft x)\triangleleft y\,,
\end{equation}
whereas the covariance (or generalized differential property) means
\begin{equation}
\label{11}
x\triangleright(ab)=(x_{(1)}\triangleright a)(x_{(2)}\triangleright b)\,,
\ \ \ 
(ab)\triangleleft x=(a\triangleleft x_{(1)})(b\triangleleft x_{(2)})\,,
\end{equation}
i.e., the $A^*$-actions respect multiplicative structure of $A$, or, in 
other words, $A$ is a left (right) $A^*$-module algebra.

One can use (e.g., left) action (\ref{9}) of $A^*$ on $A$
to define on their tensor product  $A\otimes A^*$
the cross-product algebra $A\cross A^*$ 
\cite{Sw,Ma1,SWZ2}. This is an associative 
algebra with the cross-multiplication rule given by 
\begin{equation} 
xa=(x_{(1)}\triangleright a)x_{(2)}\equiv \ <\!x_{(1)},a_{(2)}\!>
a_{(1)}x_{(2)} 
\label{12}
\end{equation}
(multiplication inside $A$ and $A^*$ does not change).
A cross-product is not a Hopf algebra but exhibits remarkable 
$A^*$-module and $A$-comodule properties~\cite{SWZ2,Sc2}. 

First, $A\cross A^*$ is covariant under the right $A^*$-action of the 
following form:  
\begin{equation} 
\mbox{\rm right $A^*$-action:} \quad\qquad
a\triangleleft x= a_{(2)}<\!x,a_{(1)}\!>\,, \ \ \ 
y\triangleleft x=y<x,1>\equiv\varepsilon (x)\,y\,,
\label{13}
\end{equation} 
to be extended on arbitrary products in $A\cross A^*$ by the covariance
condition
\begin{equation}
\label{13a}
(pq)\triangleleft x=(p\triangleleft x_{(1)})(q\triangleleft x_{(2)})\,, \ \
\ \ p,q\in (A\cross A^*)\,.
\end{equation}
Surely, this needs to be consistent with (\ref{12}). Let us check it:
$$
(ya)\triangleleft x=(y\triangleleft x_{(1)})(a\triangleleft x_{(2)})=
\varepsilon (x_{(1)})y<x_{(2)},a_{(1)}>a_{(2)}
=<x,a_{(1)}>ya_{(2)} 
$$
$$
= <x,a_{(1)}><y_{(1)},a_{(3)}>a_{(2)}y_{(2)}
=\varepsilon (x_{(2)})<x_{(1)},a_{(1)}><y_{(1)},a_{(3)}>a_{(2)}y_{(2)}
$$
\begin{equation}
\label{15}
=<y_{(1)},a_{(2)}>(a_{(1)}\triangleleft x_{(1)})(y_{(2)}
\triangleleft x_{(2)})=
(<y_{(1)},a_{(2)}>a_{(1)}y_{(2)})\triangleleft x\,.
\end{equation}

It is known~\cite{Ma1} that a covariant right action $F\triangleleft H$
of a Hopf algebra $H$ on an algebra $F$ implies a covariant left coaction
$F\rightarrow H^*\otimes F$ of the Hopf dual $H^*$ on $F$. The
correspondence is defined by
\begin{equation}
\label{14a}
f\triangleleft h=<h\,,f^{(1)}>f^{(0)}\,,
\end{equation}
where a coaction is assumed to be 
$f\rightarrow f^{(1)}\otimes f^{(0)}$ with
$h\in H, \ f^{(1)}\in H^*, \ \,f,f^{(0)}\in F$. For coactions, `covariant'
still means `respecting multiplication'. This is expressed by
\begin{equation}
\label{14b}
(fg)\longrightarrow f^{(1)}g^{(1)}\otimes f^{(0)}g^{(0)}\,.
\end{equation}
In our case, the left $A$-coaction dual to (\ref{13}) is
\begin{equation} 
\mbox{\rm left $A$-coaction:} \qquad\qquad\qquad  
a\longrightarrow\Delta(a)\equiv a_{(1)}\otimes a_{(2)}\,, \ \ \
y\longrightarrow 1\otimes y\,.
\label{14}
\end{equation}
The very last relation explains why the elements of
$A^*$ are called left-invariant in this situation. 

Further, $A\cross A^*$ is covariant under a left $A^*$-action
and also under its dual right $A$-coaction. Explicit form of the 
$A^*$-action is taken to be the well-known Hopf adjoint,
\begin{equation}
\label{18}
x\stackrel{\mbox{\scriptsize ad}}{\triangleright}p=x_{(1)}pS(x_{(2)})\,,
\ \ \ \ \ p\in A\cross A^*\,,
\end{equation}
which is evidently covariant:
\begin{equation}
\label{18a}
x\stackrel{\mbox{\scriptsize ad}}{\triangleright}(pq)
=x_{(1)}pqS(x_{(2)})=x_{(1)}pS(x_{(2)})x_{(3)}qS(x_{(4)})
=(x_{(1)}\stackrel{\mbox{\scriptsize ad}}{\triangleright}p)
(x_{(2)}\stackrel{\mbox{\scriptsize ad}}{\triangleright}q)\,.
\end{equation}
Moreover, for $p=a\in A$ one shows~\cite{Ma2,SWZ2} that
$$
x\stackrel{\mbox{\scriptsize ad}}{\triangleright}a
\equiv x_{(1)}aS(x_{(2)})
=a_{(1)}<\!x_{(1)},a_{(2)}\!>x_{(2)}S(x_{(3)})
$$
\begin{equation}
\label{17}
=\varepsilon (x_{(2)})a_{(1)}<\!x_{(1)},a_{(2)}\!>
=a_{(1)}<\!x,a_{(2)}\!>=x\triangleright a\,,
\end{equation}
i.e., we recover the left action (\ref{9}) and can rewrite (\ref{18}) as
 \begin{equation}
\mbox{\rm left $A^*$-action:} \qquad\qquad\qquad  
x\triangleright a= a_{(1)}<\!x,a_{(2)}\!>\,, \ \ \
x\triangleright y=x_{(1)}yS(x_{(2)})\,.
\label{16}
\end{equation}
The corresponding dual right $A$-coaction is deduced from the general
rule~\cite{Ma1} analogous to (\ref{14a}), which relates left action
$H\triangleright F$ with right coaction $F\rightarrow F\otimes H^*$:
\begin{equation}
\label{19a}
h\triangleright g=<h\,,g^{(1)}>g^{(0)}\,, \ \ 
g\longrightarrow g^{(0)}\otimes g^{(1)}\,, \ \  
h\in H, \ g^{(1)}\in H^*, \ \,g,g^{(0)}\in F\,,
\end{equation}
and is explicitly given by~\cite{SWZ2}
\begin{equation}
\mbox{\rm right $A$-coaction:} \qquad\qquad  
a\longrightarrow\Delta(a)\equiv a_{(1)}\otimes a_{(2)}\,, \ \ \
y\longrightarrow (e^\alpha 
\stackrel{\mbox{\scriptsize ad}}{\triangleright}y)\otimes e_\alpha \,,
\label{19}
\end{equation}
where $\{e_\alpha \},\{e^\alpha \}$ are dual bases in $A$ and $A^*$.
Note that in both (\ref{14}) and (\ref{19}) the coaction on the $A$\,-part
of $A\cross A^*$ is just a coproduct.

Being the covariant (co)actions, 
eqs. (\ref{13}),(\ref{14}),(\ref{16}) and (\ref{19}) characterize 
$A\cross A^*$ as a left (right) (co)module algebra. It is in this sense
that the cross-product algebra $A\cross A^*$ may be called 
bicovariant~\cite{SWZ2,Sc2}. Of course, this bicovariance is merely a
reflection of the underlying Hopf algebra structure of $A$.

\vspace{.3cm}

{\bf III. Woronowicz's differential complex as a Hopf algebra}

\vspace{.1cm}

Let us now recall the basic definitions of the Woronowicz noncommutative
differential calculus~\cite{Wo,AC}. First, a basis $\{\omega ^i\}$
of left-invariant 1-forms should be chosen in the bimodule $\Gamma$ in
(\ref{1}). Any element $\rho \in\Gamma$ can be uniquely represented
as $\rho =a_i\omega ^i, \ a_i\in A$. 
Next, one specifies commutation relations between functions and differential
forms,
\begin{equation}
\label{20}
\omega ^ia=(f_{j}^{i}\triangleright a)\omega ^j\,,
\end{equation}
the coalgebra structure of $\Gamma$,
\begin{equation}
\label{21}
\Delta(\omega ^i)=1\otimes \omega ^i + \omega ^j\otimes r_{j}^{i} \,,
\end{equation}
and a differential map $d : A\longrightarrow\Gamma$:
\begin{equation}
\label{22}
da=(\chi _i\triangleright a)\omega ^i \,. 
\end{equation}
Here $a$ is arbitrary element of $A$, \ $r_{j}^{i}\in A$, \ $\chi _i$
and $f_{j}^{i}$ belong to $A^*$. The Hopf-algebra consistency 
(or bicovariance) conditions of the calculus are:
\begin{eqnarray}
(\Delta\otimes id) \circ \Delta = (id\otimes \Delta) \circ \Delta
& \Longrightarrow & 
\Delta(r_{j}^{i})=r_{i}^{k}\otimes r_{k}^{j}\,, 
\label{23} \\
\omega(ab)=(\omega a)b
& \Longrightarrow & 
\Delta(f_{j}^{i})=f_{k}^{i}\otimes f_{j}^{k}\,, 
\label{24} \\
\Delta(\omega a) = \Delta(\omega)\Delta(a)
& \Longrightarrow & 
(f_{i}^{j}\triangleright a)r_{k}^{i}=r_{i}^{j}(a\triangleleft f_{k}^{i})\,, 
\label{25} \\
d(ab)=(da)b+a\,db 
& \Longrightarrow & 
\Delta(\chi_{i})=\chi_{j}\otimes f_{i}^{j}+ 1\otimes \chi_{i}\,,
\label{26} \\
\Delta \circ d = (d\otimes 1 + 1\otimes d) \circ \Delta 
& \Longrightarrow &  
a\triangleleft\chi_{i}= (\chi_{j}\triangleright a)r_{i}^{j}\,, 
\label{27}
\end{eqnarray}
supplemented by the formulas
\begin{equation}
\label{28}
\varepsilon(f_{j}^{i})=\delta_{j}^{i}\,, \ \ \ 
\varepsilon(r_{j}^{i})=\delta_{j}^{i}\,, \ \ \
S(f_{k}^{j})f_{i}^{k} = \delta _{i}^{j}\,, \ \ \ 
S(r_{i}^{k})r_{k}^{j} = \delta _{i}^{j}\,, 
\end{equation}
which are obtained from the properties of counit and antipode. Woronowicz's
theory asserts that every set of elements $\{r_{i}^{j},f_{i}^{j},\chi_i\}$
obeying eqs. (\ref{23})--(\ref{28}) gives us an example of a bicovariant 
differential calculus on the Hopf algebra $A$.

For illustration, let us derive (\ref{27}) (cf.~\cite{BGMST}). 
$$
\Delta(da)=\Delta((\chi_i\triangleright a)\omega^i)
=<\chi_i,a_{(3)}>(a_{(1)}\otimes a_{(2)})(\omega ^j\otimes r_{j}^{i}
+1\otimes \omega ^i)
$$
\begin{equation}
\label{29}
=a_{(1)}\omega ^j\otimes (\chi _i\triangleright a_{(2)})r_{j}^{i}
+a_{(1)}\otimes (\chi _i\triangleright a_{(2)})\omega ^i\,,    
\end{equation}
$$
da_{(1)}\otimes a_{(2)}+a_{(1)}\otimes da_{(2)}
=(\chi _i\triangleright a_{(1)})\omega ^i\otimes a_{(2)}
+a_{(1)}\otimes(\chi _i\triangleright a_{(2)})\omega ^i
$$
$$
=a_{(1)}<\chi _i,a_{(2)}>\omega ^i\otimes a_{(3)}
+a_{(1)}\otimes(\chi _i\triangleright a_{(2)})\omega ^i
$$
\begin{equation}
\label{30}
=a_{(1)}\omega ^i\otimes (a_{(2)}\triangleleft \chi _i)
+a_{(1)}\otimes(\chi _i\triangleright a_{(2)})\omega ^i\,.
\end{equation}
Independence of $\{\omega ^i\}$ yields
\begin{equation}
\label{31}
a_{(1)}\otimes (\chi _i\triangleright a_{(2)})r_{j}^{i}
=a_{(1)}\otimes (a_{(2)}\triangleleft \chi _j)\,.
\end{equation}
Acting on both sides of this equation by $\varepsilon \otimes id$,
we come to (\ref{27}).

Consider now the graded Hopf algebra (differential complex) 
$\Gamma^{\wedge}$ given by (\ref{1}),(\ref{3}) jointly with its dual 
$(\Gamma^{\wedge})^*$:
\begin{equation}
\label{32} 
\begin{array}{ccccccc}
   A& \stackrel{d}{\longrightarrow} &
 \Gamma& \stackrel{d}{\longrightarrow} &
 \Gamma^{2}& \stackrel{d}{\longrightarrow} &
 \ldots \\
   \updownarrow &&\updownarrow&&\updownarrow&&\\
   A^{*}& \stackrel{d^{*}}{\longleftarrow} &
 \Gamma^{*}& \stackrel{d^{*}}{\longleftarrow} &
 \Gamma^2{}^*& \stackrel{d^{*}}{\longleftarrow} &
 \ldots
 \end{array}
\end{equation}
(vertical arrows indicate non-zero duality brackets implied by grading).
Analogously to (\ref{12}), an associative algebra 
${\cal G}=\Gamma^\wedge\cross \,(\Gamma^\wedge)^*$ 
can be introduced using the cross-product construction 
(here $ (\Gamma^\wedge)^*=A^*\oplus\Gamma^*\oplus\ldots\,)$.  
We place ${\cal G}$ in the center of our approach. 
It means that we assume the following guiding principle: 

{\em All cross-commutation
relations among functions, forms, Lie derivatives, and inner 
derivations are to be chosen according to the rules (\ref{12}) of a
cross-product algebra. In other words, given 
Woronowicz's calculus (and, hence, the Hopf algebra $\Gamma^\wedge$), 
we then have to use only standard Hopf-algebra technique 
$\Gamma^\wedge\!\Rightarrow\!(\Gamma^\wedge)^*\!\Rightarrow\! 
\Gamma^\wedge\!\cross(\Gamma^\wedge)^* $ to construct the whole algebra 
of these four types of elements.} 

The resulting algebra is bicovariant by construction. Its bicovariance
in the sense of Woronowicz's left and right covariance~\cite{Wo} is implied
by the Hopf-algebra nature of $\Gamma^\wedge$~\cite{Br}, whereas its
bicovariance in the sense of Schupp, Watts and Zumino, expressed by eqs.
(\ref{13}),(\ref{14}),(\ref{16}) and (\ref{19}), proves to be an inherent
feature of the cross-product (see Sect. 2), and stems, at the very end,
from the same Hopf structure of $\Gamma^\wedge$.

\vspace{.3cm}

{\bf IV. Explicit form of commutational relations}

\vspace{.1cm}

It only remains to put all the relevant objects in the
corresponding `boxes'. We already know that functions and 1-forms are
situated in $A$ and $\Gamma$, respectively. Owing to (\ref{14}), one may 
consider $A^*$ (acting on $A$ from the left) as an algebra of 
left-invariant (and $A\cross A^*$ -- of general) vector fields on a 
quantum group $A$. It is generally accepted \cite{SWZ1,SWZ2,SWZ3,AS} 
that  Lie derivatives ${\cal L}_h$ along a (left-invariant) vector
field $h\in A^*$ must be related with its action on arbitrary elements
of ${\cal G}$:
\begin{equation}
\label{33}
{\cal L}_h\doteq h\stackrel{\mbox{\scriptsize ad}}{\triangleright} 
\end{equation} 
which, due to (\ref{17}), reduces to ordinary left action (\ref{9}) 
$h\triangleright \rho $ for $\rho \in\Gamma^\wedge$.

It seems also natural to relate inner derivations with elements of
$\Gamma^*$ \cite{Sc1}. We propose the following definition~\cite{Vl}.
Let $\gamma _i\in\Gamma^*$ be determined by fixing its duality bracket
with a general element of $\Gamma$,
\begin{equation}
<\!\gamma _i,a\omega^j\!>=\varepsilon (a)\,\delta _{i}^{j}\,, 
\label{34}
\end{equation}
and $<\!\gamma_i,\rho\!>=0$ for $\rho\in A,\Gamma^2,\Gamma^3,\ldots$\ \,.
Then we can define a basis of inner derivations
$\{\imath _i\}$ as follows:
\begin{equation}
\imath _i\doteq\gamma _i\stackrel{\mbox{\scriptsize ad}}{\triangleright} 
\label{35}
\end{equation}
(the same comment as for eq. (\ref{33}) applies). 
Here we make no attempt to associate
some $\imath _h\in\Gamma^*$ with any $h\in A^*$, for it
looks unnatural in the context of our approach (see, however, 
~\cite{SWZ3,SW} for a discussion of such a possibility).

The cross-product algebra we are seeking for, i.e., an algebra which
includes four types of differential-geometric objects,
$a, \omega ^i, {\cal L}_h$ and $\imath_i$, is implicitly contained in
the above definitions. In order to make it more transparent, we 
employ these definitions for obtaining a set of helpful relations.

To begin with, the dual differential map $d^*$ is introduced by
\begin{equation}
\label{36}
<d^*\theta ,\rho >\doteq<\theta , d\rho >\,, \ \ \ \ \
\rho \in\Gamma^\wedge\,, \ \theta \in {\Gamma^\wedge}^*\,.
\end{equation}
It commutes with elements of $A^*$,
\begin{equation}
\label{37}
d^*\circ h=h\circ d^*, \ \ \ \ \mbox{\rm i.e.,} \ \
d^*(h\theta )=h\,d^*\theta\,, \ \ \ \ h\in A^*\,,
\end{equation}
and transforms $\gamma _i$ to $\chi _i$:
\begin{equation}
\label{38}
\chi _i=d^*\gamma _i\,.
\end{equation}
Both formulas are derived via duality:
$$
<d^*(h\theta ),\rho >=<h\theta ,d\rho >=<h\otimes \theta ,\Delta(d\rho)>
=<h,\rho _{(1)}><\theta ,d\rho _{(2)}>
$$
\begin{equation}
\label{39}
=<h,\rho _{(1)}><d^*\theta ,\rho _{(2)}>
=<h\otimes d^*\theta ,\Delta(\rho )>=<h\,d^*\theta ,\rho >
\end{equation}
(we used $<h,d\rho _{(1)}>=0$), and
$$
<d^*\gamma _i,a>=<\gamma _i,da>
=<\gamma _i,(\chi _j\triangleright a)\omega ^j>
=\varepsilon (\chi _i\triangleright a)
$$
\begin{equation}
\label{40}
=\varepsilon (a_{(1)})<\chi _i,a_{(2)}>=<\chi _i,a>\,.
\end{equation}

Further, to verify that the coproduct of $\gamma _i$ is given by
\begin{equation}
\label{41}
\Delta(\gamma _i)=1\otimes \gamma _i+\gamma _j\otimes f_{i}^{j}
\end{equation}
it suffices to compute its bracket with a general element in
$ A\otimes \Gamma + \Gamma\otimes A$:
$$
<\Delta(\gamma _i)-1\otimes \gamma _i-\gamma _j\otimes f_{i}^{j}\,,\,
a\otimes b\omega ^k+c\omega ^k\otimes e>
$$
$$
=<\gamma _i,ab\omega ^k+c\omega ^ke>
-\varepsilon (a)<\gamma _i,b\omega ^k>
-<\gamma _j,c\omega ^k><f_{i}^{j},e>
$$
\begin{equation}
\label{42}
=<\gamma _i,c\,(f_{m}^{k}\triangleright e)\,\omega ^m>
-\varepsilon (c)<f_{i}^{k},e>
=\varepsilon (c)\,\varepsilon (f_{i}^{k}\triangleright e)
-\varepsilon (c)<f_{i}^{k},e>=0\,,
\end{equation}
where $a,b,c,e\in A$. A comparison of (\ref{41}) with (\ref{21}) displays
a `left appearance' of $\Delta(\gamma _i)$. Nevertheless, unlike the
$\omega ^i$\,-\,case (\ref{21}), we prefer not to use the words `left 
invariance' here, to avoid confusion with the left invariance under $A$\,-
coaction
(\ref{14}) appropriate to any object in $\Gamma^\wedge{}^*$. However, a
similarity of (\ref{41}) and (\ref{21}) enables one to show in a way quite
analogous to~\cite{Wo} that any element $\theta \in\Gamma^*$ is uniquely
represented in the form $\theta =h^i\gamma _i, \ \,h^i\in A^*$.

Now we are in a position to derive 10 commutation relations among
$a\in A, \ \omega ^i\in \Gamma$, $\chi _i\in A^*$ 
and $\gamma _i\in \Gamma^*$.
Three of them are already present in the original Woronowicz theory.
They are: internal multiplication rule inside the algebra $A$, 
eq. (\ref{20}), and the recipe
how to (wedge) multiply $\omega^i$. The latter is unambiguously
fixed in the framework of Woronowicz's scheme~\cite{Wo} but generally 
cannot be written down in a closed form (see~\cite{CWSWW,AC,Is}). Another 
four, 
\begin{equation}
\label{43}
\chi_{i}a-a\chi_{i} = (\chi_{j}\triangleright a) f_{i}^{j}\,,
\end{equation}
\begin{equation}
\label{44}
\gamma_{i}a-a\gamma_{i}=0\,,
\end{equation}
\begin{equation}
\label{45}
\chi_{i}\omega^{j}-\omega^{j}\chi_{i} = 
C_{lk}^{j}\omega^{l}f_{i}^{k}\,, \ \ \ \ 
C_{lk}^{j}\doteq<\!\chi_{k}, r_{l}^{j}\!>\,,
\end{equation}
\begin{equation}
\label{46}
\gamma_{i}\omega^{j}+\omega^{j}\gamma_{i} = f_{i}^{j}\,,
\end{equation}
are immediately obtained
by applying the cross-product rule (\ref{12}) to
$\Gamma^\wedge\!\cross(\Gamma^\wedge)^*$.
 The remaining commutation relations require the use 
of the duality arguments. Let us first derive a formula 
\begin{equation} 
\label{47} 
\gamma_{i}h = (r_{i}^{j} \triangleright h )\gamma _j \equiv 
<h_{(2)},r_{i}^{j}>h_{(1)}\gamma _j\,.  
\end{equation} 
We have 
$$ 
<\gamma_{i}h\,,a\omega ^k>
=<\gamma_{i}\otimes h\,,a_{(1)}\otimes a_{(2)}\omega ^k
+a_{(1)}\omega ^j\otimes a_{(2)}r_{j}^{k}>
$$
$$
=<\gamma_{i}\,,a_{(1)}\omega ^j><h\,,a_{(2)}r_{j}^{k}>
=\delta_{i}^{j}\varepsilon(a_{(1)})<h_{(1)},a_{(2)}><h_{(2)},r_{j}^{k}>
=<r_{i}^{k}\triangleright h\,,a>
$$
\begin{equation}
\label{48}
=<r_{i}^{j}\triangleright h\,,a_{(1)}>\varepsilon(a_{(2)})\delta_{j}^{k} 
=<r_{i}^{j}\triangleright h\,,a_{(1)}><\gamma_j,a_{(2)}\omega ^k> 
=<(r_{i}^{j}\triangleright h)\gamma_j\,,a\omega ^k>\,. 
\end{equation}
Using (\ref{37}) and (\ref{38}), we come to analogous formula for $\chi_i$:
\begin{equation}
\label{49}
\chi_{i}h = (r_{i}^{j}\triangleright h)\chi _j
\equiv <h_{(2)},r_{i}^{j}>h_{(1)}\chi_j\,.
\end{equation}
This can be also proved by a direct calculation:
$$
<\chi _ih\,, a>=<\chi _i,a_{(1)}><h\,,a_{(2)}>
=<h\,,\,a\triangleleft\chi _i>
=<h\,,\,(\chi _j\triangleright a)r_{i}^{j}>
$$
$$
=<h_{(1)},\chi _j\triangleright a><h_{(2)},r_{i}^{j}>
=<h_{(1)},a_{(1)}><\chi _j,a_{(2)}><h_{(2)},r_{i}^{j}>
$$
\begin{equation}
\label{50}
=<r_{i}^{j}\triangleright h\,,a_{(1)}><\chi _j,a_{(2)}>
=<(r_{i}^{j}\triangleright h)\chi _j,a>\,. 
\end{equation}
It is worth mentioning that the same technique leads to a helpful
formula
\begin{equation}
\label{51}
f_{i}^{j}h=(r_{i}^{k}\triangleright h\triangleleft S^{-1}(r_{m}^{j}))
f_{k}^{m}\equiv \ <h_{(1)},S^{-1}(r_{m}^{j})>h_{(2)}f_{k}^{m}
<h_{(3)},r_{i}^{k}>
\end{equation}
which can be used, in conjunction with (\ref{49}), to deduce
the structure relations of bicovariant differential calculus in the form
given in~\cite{AC,Be}:
\begin{equation}
\label{52}
\chi_{i}\,\chi_{j}-\sigma_{ij}^{lk}\chi_{l}\,\chi_{k} =
C_{ij}^{k}\,\chi_{k}\,, \ \ \ \ 
\sigma _{ij}^{kl}\doteq<f_{j}^{k}\,,r_{i}^{l}>\,,
\end{equation}
\begin{equation}
\label{53}
\sigma _{ij}^{mn}f_{k}^{i}f_{l}^{j}
=\sigma _{kl}^{ij}f_{i}^{m}f_{j}^{n}\,,
\end{equation}
\begin{equation}
\label{54}
\chi _kf_{l}^{n}=\sigma _{kl}^{ij}f_{i}^{n}\chi _j\,,
\end{equation}
\begin{equation}
\label{55}
C_{mn}^{i}f_{j}^{m}f_{k}^{n}+f_{j}^{i}\chi _k
=\sigma _{jk}^{mn}\chi _mf_{n}^{i}+C_{jk}^{m}f_{m}^{i}\,.
\end{equation}

Now we can list the remaining three commutational relations.
One of them is (\ref{52}), and the other two are as follows:
\begin{equation}
\label{56}
\gamma_{i}\chi_{j}-\sigma_{ij}^{lk}\chi_{l}\gamma_{k} =
C_{ij}^{k}\gamma_{k}\,,
\end{equation}
\begin{equation}
\label{57}
<\!\gamma_{i}\gamma_{j}\,,a\,\omega^{m}\omega^{n}\!>
=\varepsilon(a)(\sigma_{ij}^{mn}-\delta_{i}^{m}\delta_{j}^{n})\,.
\end{equation}
Eq. (\ref{56}) stems from (\ref{47}), whereas (\ref{57})
is verified by a straightforward calculation.

Thus, we have completed the explicit construction of the cross-product 
algebra generated by $a, \omega ^i, \chi _i$ and $\gamma _i$.

\newpage

{\bf V. Cartan identity}

\vspace{.1cm}

Remarkably, this quantum algebra exhibits some features which exactly
correspond to the well-known classical relations. First, the Lie
derivatives commute with exterior differentiation:
\begin{equation}
{\cal L}_h\circ d=d\circ{\cal L}_h \,, \ \ \ \ h\in A^*\,. 
\label{58}
\end{equation}
Really, prove it for $a\in A$\,:
\begin{equation}
\label{59}
h\triangleright da=da_{(1)}<h\,,a_{(2)}>=d(h\triangleright a)\,.
\end{equation}
Then, from
\begin{equation}
\label{60}
h\triangleright(db_1\ldots db_n)=(h_{(1)}\triangleright 
db_1)\ldots (h_{(n)}\triangleright 
db_n)=d(h_{(1)}\triangleright b_1)\ldots 
d(h_{(n)}\triangleright b_n) 
\end{equation}
and the Leibniz rule it follows that 
\begin{equation}
\label{61}
h\triangleright(d(a\,db_1\ldots db_n))= 
d(h\triangleright(a\,db_1\ldots db_n))\,,  
\end{equation}
which is exactly (\ref{58}).

Furthermore, the Cartan identity in the classical form
can be shown to be valid:
\begin{equation}
\label{62}
{\cal L}_{\chi _i}=d\circ\imath _i+\imath _i\circ d\,. 
\end{equation}
One needs to verify that
\begin{equation}
\chi _i\triangleright\rho =d(\gamma _i\triangleright\rho )+
\gamma _i\triangleright(d\rho )\,, \ \ \ \ \rho \in\Gamma^\wedge\,.
\label{63}
\end{equation}
For $\rho =a\!\in\! A$ eq. (\ref{63}) is almost trivial and follows from 
\begin{equation}
\label{64}
\gamma _i\triangleright a=0\,, \ \ \ \ \
\gamma _i\triangleright da=\chi_i\triangleright a\,. 
\end{equation}
Let $\rho =a\,db \ \ (a,b\in A)\,.$ To show that
\begin{equation}
\chi _i\triangleright(a\,db)=d(\gamma _i\triangleright a\,db)+
\gamma _i\triangleright(da\,db)
\label{65}
\end{equation}
we calculate each term separately,
$$
\chi _i\triangleright(a\,db)=a_{(1)}db_{(1)}<\!1\otimes \chi _i+
\chi _j\otimes f_{i}^{j}\,,\ a_{(2)}\otimes b_{(2)}\!> 
$$
\begin{equation}
\label{66}
=a(\chi _i\triangleright db)+(\chi _j\triangleright a)(f_{i}^{j}
\triangleright db)\,,
\end{equation}
\begin{equation}
\label{67}
d(\gamma _i\triangleright(a\,db))=da\,(\chi _i\triangleright b)+
a\,d(\chi _i\triangleright b)\,, 
\end{equation}
\begin{equation}
\label{68}
\gamma _i\triangleright(da\,db)=-da\,(\chi _i\triangleright b)+
(\chi _j\triangleright a)(f_{i}^{j}\triangleright db)\,, 
\end{equation}
and then use (\ref{59}).

At last, consider the general case $\rho =a\,db\,B\,, \ 
B=dc_1\ldots dc_n\,,
\ a,b,\ldots ,c_i\in A\,$:
\begin{equation}
\label{69}
\chi _i\triangleright(a\,db\,B)=a\,db\,(\chi _i\triangleright B)+
a(\chi _j\triangleright db)(f_{i}^{j}\triangleright B) 
 +(\chi _k\triangleright a)(f_{j}^{k}\triangleright db)
(f_{i}^{j}\triangleright B)\,, 
\end{equation}
\begin{equation}
\label{70}
\gamma _i\triangleright(a\,db\,B)=-a\,db\,(\gamma _i\triangleright 
B)+ a(\chi _j\triangleright b)(f_{i}^{j}\triangleright B)\,, 
\end{equation}
$$
d(\gamma _i\triangleright(a\,db\,B))=-da\,db\,(\gamma _i\triangleright 
B)+ a(\chi _j\triangleright db)(f_{i}^{j}\triangleright B) 
$$
\begin{equation}
\label{71} 
+a\,db\,d(\gamma _i\triangleright B) + da\,(\chi _j\triangleright 
b)(f_{i}^{j}\triangleright B) + a(\chi _j\triangleright 
b)\,d(f_{i}^{j}\triangleright B)\,, 
\end{equation}
$$ 
\gamma_i\triangleright(da\,db\,B)=da\,db\,(\gamma _i\triangleright B) 
-da\,(\chi _j\triangleright b)(f_{i}^{j}\triangleright B) 
$$ 
\begin{equation}
\label{72}
+(\chi_k\triangleright a)(f_{j}^{k}\triangleright db) 
(f_{i}^{j}\triangleright B)\,. 
\end{equation}
After summing this up, it remains to prove that
\begin{equation}
\label{73}
a\,db\,(\chi _i\triangleright B)=a\,db\,d(\gamma _i\triangleright 
B)\,, 
\end{equation}
or
\begin{equation}
\label{74}
\chi _i\triangleright(dc_1\ldots dc_n)=d(\gamma _i\triangleright
(dc_1\ldots dc_n))\,, 
\end{equation}
that is the same problem at a lower level. Thus, the proof is completed 
by induction.

To conclude this section, we compare the duality
$<\Gamma^*,\Gamma>$ used above (`vertical' duality in (\ref{32}) between
1-forms and inner derivations) with a duality 
$<\!<\!A^*,\Gamma\!>\!>$ between 
vector fields $\in\!A^*$ and differential 1-forms $\in\!\Gamma$.
The latter is a natural ge\-ne\-ralization of ordinary classical duality
and is assumed as a basis of an alternative construction of bicovariant 
differential calculus on the Hopf algebras in~\cite{AS}. It is easily 
seen that the dual differential map $d^*$ establishes a direct relation 
between these two dualities in the following way:  
\begin{equation} 
\label{75}
<\!<d^*\theta  ,\rho>\!> \ = \ <\theta ,\rho >\,, \ \ \ \
\theta \in\Gamma^*, \ \,\rho \in\Gamma\,.  
\end{equation}

\vspace{.3cm}

{\bf VI. Comparison with other approaches}

\vspace{.1cm}

Now the above results (mostly, the commutation relations 
(\ref{43})--(\ref{46}), (\ref{52}), (\ref{56}) and (\ref{57}))
are to be compared with other approaches known in the literature
~\cite{SWZ1,SWZ3,Sc1,Is}.
To achieve this, it is convenient to chose another set of generators
for the $\Gamma^\wedge{}^*$-part of
our cross-product algebra. We switch from $\chi _i,\gamma _i$ to
$\tilde\chi _i,\tilde\gamma _i$ defined by
\begin{equation}
\label{76}
da = \omega^{i}(\tilde\chi_{i} \triangleright a)\,,
\end{equation}
\begin{equation}
\label{77}
<\! \tilde\gamma_{i}, \omega^{j}a \!> = \varepsilon(a)\delta^{j}_{i}\,, \ \
\ <\! \tilde\gamma_{i}, \rho \!>=0\,, \ \ \ \ \rho\in A,\Gamma^{2},
 \Gamma^{3},\ \cdots \ \ .
\end{equation}
Introducing also $\varphi _{j}^{i}\in A^*$ via
\begin{equation}
\label{79}
a\omega ^j=\omega ^i(\varphi _{i}^{j}\triangleright a)
\end{equation}
and proceeding by complete analogy with Sect. 3 and 4, we obtain
\begin{equation}
\label{80}
\varphi_{i}^{j} = S^{-1}(f_{i}^{j})\,, \ \ \   
\tilde\gamma_{i} = \sigma_{mk}^{mj}\gamma_{j}\varphi_{i}^{k}\,, \ \ \ 
\tilde\chi_{i} = \varphi_{i}^{j}\chi_{j}=d^*\tilde\gamma _i\,,
\end{equation}
\begin{equation}
\label{81}
\Delta(\varphi _{j}^{i})=\varphi _{j}^{k}\otimes \varphi _{k}^{i}\,, \ \ \
\Delta(\tilde\gamma_{i}) = \tilde\gamma_{i}\otimes 1 
+ \varphi_{i}^{j}\otimes \tilde\gamma_{j}\,, \ \ \ 
\Delta(\tilde\chi_{i}) = \tilde\chi_{i}\otimes 1 
+ \varphi_{i}^{j}\otimes \tilde\chi_{j}\,.
\end{equation}

As for commutational relations, in the 
$\{a,\omega\}$\,-\,sector they remain 
unchanged, those between $a,\omega$ and $\chi,\gamma$ follow directly from 
(\ref{12}), 
\begin{equation}
\label{82}
 \tilde\chi_{i}a - (\varphi_{i}^{j} \triangleright a) 
\tilde\chi_{j} = \tilde\chi_{i} \triangleright a\,,
\end{equation}
\begin{equation}
\label{83}
\tilde\gamma_{i}a-(\varphi_{i}^{j}\triangleright a)\tilde\gamma_{j} = 0\,, 
\end{equation}
\begin{equation}
\label{84}
\tilde\chi_{i}\,\omega^{j}-\tilde\sigma_{ik}^{jl}\,\omega^{k}\tilde\chi_{l} 
=\tilde C_{ki}^{j}\,\omega^{k}\,,
\end{equation} 
\begin{equation}
\label{85}
\tilde\gamma_{i}\,\omega^{j} + 
\tilde\sigma_{ik}^{jl}\,\omega^{k}\tilde\gamma_{l} = \delta_{i}^{j}\,,
\end{equation}
where 
\begin{equation}
\label{86}
\tilde\sigma_{jl}^{ki} \doteq <\! \varphi_{j}^{i}, r_{l}^{k} \!> = 
(\sigma^{-1})_{jl}^{ki}\,, \ \ \
\tilde C_{jk}^{i} \doteq <\! \tilde\chi_{k}, 
r_{j}^{i} \!> = C_{sl}^{i}(\sigma^{-1})_{kj}^{sl}\,,
\end{equation}
and those inside $\{\chi ,\gamma \}$ look like
\begin{equation}
\label{87}
\tilde\gamma_{i}\tilde\chi_{j} - 
\tilde\sigma_{ji}^{kl}\tilde\chi_{l}\tilde\gamma_{k} 
= \tilde C_{ij}^{k}\tilde\gamma_{k}\,,
\end{equation}
\begin{equation}
\label{88}
\tilde\chi_{i}\tilde\chi_{j} - 
\tilde\sigma_{ji}^{kl}\tilde\chi_{l}\tilde\chi_{k} 
= \tilde C_{ij}^{k}\tilde\chi_{k}\,,
\end{equation}
\begin{equation}
\label{89}
<\!\tilde\gamma_{i}\tilde\gamma_{j}\,, 
\omega^{m}\omega^{n}a\!>=\varepsilon(a)(-\tilde\sigma_{ji}^{mn}+ 
\delta_{j}^{m}\delta_{i}^{n})\,.
\end{equation}
Formulas (\ref{87}),(\ref{88}) are obtained with the use of
\begin{equation}
\label{90}
\tilde\chi _i\,h=(h\triangleleft r_{i}^{j})\tilde\chi _j\,, \ \ \ \
\tilde\gamma _i\,h=(h\triangleleft r_{i}^{j})\tilde\gamma _j\,,
\end{equation}
that can be derived similarly to (\ref{47}),(\ref{49}).

The resulting cross-commutation algebra conforms to Schupp's 
paper~\cite{Sc1}.

\vspace{.3cm}

{\bf VII. $R$-matrix formulation of differential calculus on $GL_q(N)$}

\vspace{.1cm}

To compare our formulas with analogous relations in~\cite{SWZ1},
we consider a specific rea\-li\-zation~\cite{Ju} of Woronowicz's 
differential calculus in case of the quantum group $GL_q(N)$, and use 
the matrix representations for all generators. Here $A,A^*$ will be the 
dual Hopf algebras~\cite{FRT} described by the relations 
\begin{equation}
\label{91}
R_{12}T_1T_2=T_2T_1R_{12}\,,\ \ \ \Delta (T)=T\otimes T\,,
\ \ \ \varepsilon (T)={\bf 1}\,,
\end{equation}
\begin{equation}
\label{92}
R_{12}L^{\pm}_2L^{\pm}_1  =  L^{\pm}_1L^{\pm}_2R_{12}\,, \ \ \  
R_{12}L^+_2L^-_1  = L_1^-L^+_2R_{12}\,,  
\end{equation}
\begin{equation}
\label{93}
\Delta(L^{\pm})  =  L^{\pm}\otimes L^{\pm} , \ \ \
\varepsilon(L^{\pm}) = {\bf 1}\,,
\end{equation}
\begin{equation}
\label{94}
<T_1,L_2^+>=R_{12}\,, \ \ \ <T_1,L_2^->=R_{21}^{-1}\,,
\end{equation}
\begin{equation}
\label{94a}
<T_1,S(L_2^+)>=R_{12}^{-1}\,, \ \ \ <T_1,S(L_2^-)>=R_{21}
\end{equation}
(generators $t_i^j\in A$ and ${l^{\pm}}_{i}^{j}\in A^*$ 
form matrices $T$ and $L^{\pm}$,
respectively), where $R$ is a special numerical matrix related to
$GL_q(N)$~\cite{FRT} which obeys the Yang-Baxter 
\begin{equation}
\label{95}
R_{12}R_{13}R_{23}=R_{23}R_{13}R_{12}
\end{equation}
and Hecke
\begin{equation}
\label{96}
R{}_{pq}^{ij}=(R^{-1})_{qp}^{ji}+\lambda\delta _{q}^{i}\delta _{p}^{j}
\ \ \ \ \ \ \ \ (\lambda =q-q^{-1})
\end{equation}
conditions. Let us also introduce a numerical matrix $D$ by
\begin{equation}
\label{97}
D_{j}^{i}\doteq\tilde R{}_{jm}^{mi}\,, \ \ \ 
R{}_{pn}^{mj}\tilde R{}_{mq}^{in}=\tilde R{}_{pn}^{mj}R{}_{mq}^{in}
=\delta _{p}^{i}\delta _{q}^{j}\,,
\end{equation}
and fix the differential map $d : A\rightarrow\Gamma$ via
\begin{equation}
\label{98}
d\,T=T\,\Omega
\end{equation}
in terms of left-invariant Maurer-Cartan forms $\Omega$.
Then the Woronowicz bicovariant differential calculus on $GL_q(N)$
is produced by the following choice~\cite{Ju,AC,Is} of the
elements $r,f,\chi$\,:
\begin{equation}
\label{102}
r_{kj}^{li}=S(t_{k}^{i})t_{j}^{l}\,, \ \ \
f_{jk}^{il}={l^-}_{k}^{i}S({l^+}_{j}^{l})\,, \ \ \ 
\chi _{k}^{l}
=\frac{1}{\lambda}\,[\,(D^{-1})_{k}^{l}-(D^{-1})_{i}^{j}f_{jk}^{il}\,]\,,
\end{equation}
which serve to define the Hopf and differential structure of the calculus
as follows (note doubling the indices due to the matrix format used):
\begin{equation}
\label{99}
\Delta(\Omega_{j}^{i})=1\otimes \Omega_{j}^{i}
+\Omega_{l}^{k}\otimes r_{kj}^{li}\,,
\end{equation}
\begin{equation}
\label{100}
\Omega_{j}^{i}t_{n}^{m}
=(f_{jk}^{il}\triangleright t_{n}^{m})\Omega_{l}^{k}\,,
\end{equation}
\begin{equation}
\label{101}
dt_{n}^{m}=(\chi _{q}^{p}\triangleright t_{n}^{m})\Omega_{p}^{q}
=t_{k}^{m}\Omega_{n}^{k}
\end{equation}
(the last equation implies 
$<\chi _{q}^{p},t_{n}^{m}>=\delta _{n}^{p}\delta _{q}^{m}$).
From (\ref{80}) and (\ref{102}) we get
\begin{equation}
\label{103}
\varphi _{jk}^{il}=l^+{}_{j}^{l}\,S^{-1}({l^-}_{k}^{i})\,, \ \ \
\tilde\chi _{k}^{l}=
\frac{1}{\lambda}\,[\,(D^{-1})_{j}^{i}\varphi_{ik}^{jl}-(D^{-1})_{k}^{l}\,]
\end{equation}
and can now write down all commutational relations in the matrix form.

If, before doing so, we perform one more redefinition,
\begin{equation}
\label{104}
J=-\tilde\gamma D\,, \ \ \
X=-\tilde\chi D\,, \ \ \
Y=1-\lambda X\,,
\end{equation}
so that     $ Y_{j}^{i}={l^+}_{k}^{i}S({l^-}_{j}^{k}) $,
we end up with a complete set of commutation relations in terms
of matrices $T,\Omega,Y$ and $J$:
\begin{equation}
\label{106}
R_{12}T_1T_2=T_2T_1R_{12}\,,
\end{equation}
\begin{equation}
\label{107}
\Omega_1T_2=T_2R_{12}^{-1}\Omega_1R_{21}^{-1}\,,
\end{equation}
\begin{equation}
\label{108}
\Omega_1R_{21}^{-1}\Omega_2R_{21}=-R_{21}^{-1}\Omega_2R_{12}^{-1}\Omega_1\,,
\end{equation}
\begin{equation}
\label{109}
Y_1T_2=T_2R_{21}Y_1R_{12}\,,
\end{equation}
\begin{equation}
\label{110}
\Omega_1R_{12}Y_2R_{21}=R_{12}Y_2R_{21}\Omega_1\,,
\end{equation}
\begin{equation}
\label{111}
J_1T_2=T_2R_{21}J_1R_{12}\,,
\end{equation}
\begin{equation}
\label{112}
\Omega_1R_{12}J_2R_{21}+R_{12}J_2R_{21}\Omega_1
=\frac{1}{\lambda}(1-R_{12}R_{21})\,,
\end{equation}
\begin{equation}
\label{113}
Y_1R_{12}Y_2R_{21}=R_{12}Y_2R_{21}Y_1\,,
\end{equation}
\begin{equation}
\label{114}
J_1R_{12}Y_2R_{21}=R_{12}Y_2R_{21}J_1\,,
\end{equation}
\begin{equation}
\label{115}
J_1R_{12}J_2R_{21}=-R_{21}^{-1}J_2R_{21}J_1\,.
\end{equation}

Several comments are in order. In this specific realization of the
Woronowicz calculus, it proves possible to present multiplication
relations for $\Omega$ in a closed form (\ref{108}). The commutation rule
(\ref{113}) for $Y$ is often called the {\em reflection 
equation}~\cite{KSk,KSa,Ma3}, and the related formula for $X$
\begin{equation}
\label{116}
X_1R_{12}X_2R_{21}-R_{12}X_2R_{21}X_1
=\lambda^{-1}(X_1R_{12}R_{21}-R_{12}R_{21}X_1)
\end{equation}
-- the {\em quantum Lie algebra}~\cite{Wo,Ma4,SWZ2,DH,LS}, because it 
generalizes classical commutator in the Lie algebra of left-invariant
vector fields. In terms of $T$ and $Y$, the left and right $A$-coactions
in (\ref{14}) and (\ref{19}) take the form
\begin{equation}
\label{117}
\mbox{\rm left:}\qquad\qquad\qquad\qquad
t_{j}^{i}\longrightarrow t_{k}^{i}\otimes t_{j}^{k}\,, \ \ \
Y_{j}^{i}\longrightarrow {\bf 1}\otimes Y_{j}^{i}\,,
\end{equation}
\begin{equation}
\label{118}
\mbox{\rm right:}\qquad\qquad\qquad
t_{j}^{i}\longrightarrow t_{k}^{i}\otimes t_{j}^{k}\,, \ \ \
Y_{j}^{i}\longrightarrow Y_{l}^{k}\otimes S(t_{k}^{i})t_{j}^{l}\,.
\end{equation}
This shows explicitly that algebra $A^*$  is
left-invariant and right-coadjoint-covariant.

Algebra (\ref{106})--(\ref{115}) is exactly the $GL_q(N)$ bicovariant 
differential algebra found in~\cite{SWZ1} and discussed further 
in~\cite{Is}. We have shown that it is produced just by application of the
cross-product recipe to the original Woronowicz differential complex,
whose Hopf-algebra properties account for bicovariance of the algebra.

\vspace{.3cm}

{\bf Acknowledgements}

\vspace{.1cm}

We are grateful to G.\,Arutyunov, T.\,Brzezinski, \v{C}.\,Burdik, 
R.\,Giachetti, A.\,Isaev, 
P.\,Pyatov, J.\,Riembelinski and M.\,Tarlini for discussions. 

This work was supported in part 
by Russian Foundation for Basic Research.

\end{document}